# Co-Optimization Scheme for Distributed Energy Resource Planning in Community Microgrids

C. Yuan, *Member, IEEE*, M. S. Illindala, *Senior Member, IEEE, and* A. S. Khalsa, *Member, IEEE*

*Abstract*—**Microgrids with distributed energy resources are being favored in various communities to lower the dependence on utility-supplied energy and cut the $CO_2$ emissions from coal-based power plants. This paper presents a co-optimization strategy for distributed energy resource planning to minimize total annualized cost at the maximal fuel savings. Furthermore, the proposed scheme aids the community microgrids in satisfying the requirements of U.S. Department of Energy (DOE) and state renewable energy mandates. The method of Lagrange multipliers is employed to maximize fuel savings by satisfying Karush-Kuhn-Tucker conditions. With the Fourier transform and particle swarm optimization, the right mix of distributed energy resources is determined to decrease the annualized cost. A case study to test the proposed scheme for a community microgrid is presented. To validate its effectiveness, an economic justification of the solution and its comparison with HOMER Pro are also illustrated.**

## I. NOMENCLATURE

*Indices:*

| | |
|---|---|
| $e$ | Index for electrical energy |
| $i, k$ | Index for time point |
| $l$ | Index for load profile |
| $nl$ | Index for net load |
| $nom$ | Index for nominal power rating |
| $re$ | Index for renewable energy |
| $st$ | Index for stochastic value |
| $t$ | Index for time |
| $th$ | Index for thermal energy |
| $thr$ | Index for threshold power rating |

*Parameters:*

| | |
|---|---|
| $\alpha$ | Coefficient of annualized $CO_2$ emissions reduction in ton/MW |
| $\beta$ | Coefficient of annualized system energy savings in MWh/MW |
| $\gamma$ | Coefficient of annualized fuel savings in MMBtu/MW |
| $\theta$ | Power capacity ratio of photovoltaic (PV) |
| $\rho$ | Power capacity ratio of renewable energy resources |
| $\omega$ | Ratio of $CO_2$ emissions reduction |
| $\delta$ | Ratio of system energy efficiency increase |
| $\mu$ | Mean value of the normal distribution |
| $\sigma^2$ | Variance of the normal distribution |
| $a, b$ | Shape parameters of Beta distribution function |
| $\tau$ | Scale parameter of Weibull distribution function |
| $K$ | Shape parameter of Weibull distribution function |
| $E_{CO_2}$ | Base amount of annual $CO_2$ emissions in tons |
| $E_l$ | Load electrical energy annual demand in MWh |
| $E_{th}$ | Load thermal energy annual demand in MWh |
| $F$ | Fuel price in $/MWh |
| $L$ | Load demand average in MW |
| $N_S$ | Total sampling points in a day |
| $P_R$ | Power rating |
| $R$ | Reserve portion of combined heat and power (CHP) |
| $T$ | Sampling time in hour |
| $V_{in}$ | Cut-in wind speed |
| $V_{out}$ | Cut-out wind speed |
| $V_R$ | Rated wind speed |
| $\eta$ | Efficiency of energy conversion |
| $r$ | Discount rate |
| $y$ | Years of lifetime |

*Variables:*

| | |
|---|---|
| $ir$ | Solar irradiance |
| $v$ | Wind speed |

C. Yuan is with GEIRI North America, San Jose, CA 95134, USA (e-mail: chen.yuan@geirina.net).

M. S. Illindala is with The Ohio State University, Columbus, OH 43210, USA (e-mail: millindala@ieee.org).

A. S. Khalsa is with the Dolan Technology Center, American Electric Power, Groveport, OH 43125, USA (e-mail: askhalsa@aep.com).

| $P$ | Distributed energy and storage resource power in MW |
|---|---|
| $E$ | Battery energy storage system (BESS) energy in MWh |
| $SOC$ | BESS state of change |
| $C$ | Annual cost in $ |
| $ER$ | $CO_2$ emissions reduction in tons |
| $SES$ | System energy savings in MWh |
| $FS$ | Fuel savings in MMBtu |
| $URe$ | Utilization of renewable energy |

## II. INTRODUCTION

IN the United States, electricity generation from fossil fuel plants has been the biggest cause of air pollution emissions. Fossil fuels such as coal and natural gas account for 66.9% of electricity generation [1]. According to the United States Environment Protection Agency (EPA), around 31% of total greenhouse gas emissions and 37% of total $CO_2$ emissions came from the electricity sector in 2014 [2]. Within the electricity sector, coal-based power plants produce 77% of $CO_2$ emissions while generating only 38.5% of the total electricity [1]. Though the problems are alleviated to a certain extent when the coal power plants approach their retirement, it is important to further cut the emissions and increase energy efficiency. Therefore, a movement toward renewable energy resources is gaining more attention. To achieve this goal, emission-reduction targets have been established [3]–[5].

The transmission and distribution (T&D) losses in power delivery to residential, commercial, and industrial consumers account for 6% of gross generation [6]. This makes a strong case for installing distributed energy resources (DERs) near the load centers. Microgrids formed by integrating DERs can help in improving the fuel savings, energy efficiency, reliability and resiliency [7].

A review of hybrid renewable energy power generation systems affecting the energy sustainability was presented in [8]. It also discussed important issues and challenges in the system design and energy management. References [9], [10] focused on minimizing the costs and $CO_2$ emissions by solving multi-objective unit commitment problem. However, the costs and emissions are not measured on a similar scale, and hence it is difficult to combine them effectively. With the penetration of intermittent renewables, the distributed power generation is becoming uncertain and depending on environmental factors. As the power balance is affected within the distribution network, battery energy storage offers to smooth the tie-line power fluctuations and compensate the power mismatch [11]–[14]. A comprehensive system planning with the optimal mix of DERs and energy storage would guarantee the power balance between demand and supply.

A multi-objective genetic algorithm was presented in [15] to optimize a dc microgrid considering capital cost and intermittency of renewables. Besides, they analyzed power supply availability vs. cost and examined the influence of sampling time. However, the authors did not exclude conditions of negative power output from renewable energy resources. In addition, a normal distribution was assumed for renewables and load demand, which is not an accurate representation. Etezadi-Amoli *et al.* used a stochastic framework to optimally size a hybrid system for minimizing the system cost and satisfying reliability requirement [16]. Reference [17] introduced a Fourier transform based sizing scheme for energy storage and diesel generators within an isolated microgrid. For minimizing cost, they recommended using gradient search method to pick the suitable mix of energy storage and diesel generators. However, the cost function is non-convex, and hence the gradient search is not well suited to obtain the global optimum point. In addition, the net load profile did not consider any stochastic information and instead use the data sampling of a single day. A genetic algorithm based optimization for DERs allocation was presented in [18] to maximize cost savings in system upgrades deferral, including the effect of annual energy loss and load interruption. Reference [19] proposed an adaptive genetic algorithm based DER allocation strategy that minimizes network power loss and node voltage deviation. Uncertainties of load and generation are modeled using fuzzy-based approach. However, references [18], [19] combined the various problems into a single multi-objective problem with different weights for each sub-objective. A major challenge of such weighted multi-objective method is in the selection of suitable criteria to determine the weights. Hence, an acceptable solution might not be guaranteed by selecting weights before conducting an optimization process.

For the energy storage, reference [20] presented the battery control strategies and sizing algorithms. Makarov *et al.* used discrete-time Fourier transform (DTFT) to decompose the power imbalance into four cyclic components for optimal use of the battery energy storage system (BESS) [21]. However, no optimization was carried out for determining the cut-off frequencies in sizing the BESS. Moreover, the stochastic characteristic of renewable energy and load demand was not considered in [20], [21]. In addition, potential investments in BESS should be carefully analyzed for achieving operational profits. The optimal sizing of energy storage for minimizing islanded microgrid operation costs and achieving maximum benefits in grid-connected operation was presented in [22]. It uses a fixed step size to compute all costs within feasible region. However, this can be time-consuming if the search region is large or the step size is tiny. The accuracy of results is affected if the step size is increased.

In summary, most of the previous works focused on cost minimization and benefit maximization. They only considered sizing certain types of DERs under external deterministic conditions. Few, if any, investigated the influence of environmental constraints and government renewable mandates. By contrast, this paper presents an algorithm to size all the diverse kinds of DERs in a community microgrid while satisfying the various regulatory constraints.

This paper presents a co-optimization scheme to size the capacities of DERs for community microgrids while meeting the U.S. Department of Energy (DOE) requirements and state renewable energy mandates. Loads and renewables are modeled based on historical data and their stochastic characteristics. The proposed scheme sizes DERs and BESS in three steps. Lagrange multipliers and DTFT are employed to conduct the sizing process. Particle swarm optimization (PSO) is used for finding the cut-off frequency and achieving global minimum cost. A case study is presented to demonstrate the benefits for critical loads in a community. Eight representative days are selected based on historical data. Evaluation indices for renewable energy utilization are presented for assessing the performance. Besides, the impact of regulatory constraints is investigated, which has not been explored earlier. Furthermore, the proposed scheme is compared and validated with results obtained using HOMER Pro.

The rest of the paper is organized as follows. Section III presents the co-optimization scheme for sizing distributed energy resources. Stochastic models and net load are also elaborated in this section. A case study is shown in Section IV to offer uninterruptible power support to critical loads and demonstrate the benefits of a community microgrid. Besides, the impact of regulatory constraints is investigated for economic justification. In Section V, the proposed DER planning scheme is compared against Homer Pro. Finally, the conclusion is given in Section VI.

## III. CO-OPTIMIZATION FOR DISTRIBUTED ENERGY RESOURCE PLANNING

This section presents a co-optimization scheme for sizing of DERs in a microgrid. The various types of DERs are determined by using methods of qualitative functional deployment (QFD) [23] and levelized cost of energy (LCOE) with the consideration of environmental factors and fuel consumptions [24], [25]. The LCOE comparison is presented in Fig. 1. As seen in Fig. 1, with the renewable electricity tax credit [26], photovoltaic and wind turbine have lower LCOE than fuel-based resources when the annual capacity factor is larger than 0.2, widely applicable in many regions. Thus, the selected DERs include renewables like solar/photovoltaic (PV) systems, wind turbines (WT), biomass powered combined heat and power (CHP) units, natural gas (NG) based CHPs and BESS.

Fig. 2 shows the flowchart of the proposed scheme to meet the regulations while minimizing the total cost with fuel savings. In step 1, the renewable resources are main variables and sized first. The threshold power rating of natural gas CHPs, i.e., $\Sigma P_{NG}^{thr}$, is determined as a reference value. Then a net load profile is formulated based on the historical information and stochastic characteristics of renewables and load demand. In step 2, the natural gas CHPs are sized together with BESS to share the net load and minimize annualized cost. The methods of DTFT and PSO are used to find optimal capacities for CHPs and BESS. In the case that the calculated total power capacity of natural gas CHPs is much larger than the threshold value, i.e. $\Sigma P_{NG} \geq \Sigma P_{NG}^{thr}$, resulting in the violation of renewable mandates constraints, the suboptimal cut-off frequency is chosen to decrease the CHP capacity. Then the procedure for sizing of CHP and BESS is repeated. However, if $\Sigma P_{NG} < \Sigma P_{NG}^{thr}$, the suboptimal cut-off frequency is chosen to increase the CHP capacity, and the procedure for sizing CHP and BESS is repeated. Only when both $\Sigma P_{NG} \geq \Sigma P_{NG}^{thr}$ and the renewable mandates are satisfied, the optimal solution is reached. This is covered in step 3 as parity check. A detailed explanation is presented as follows.

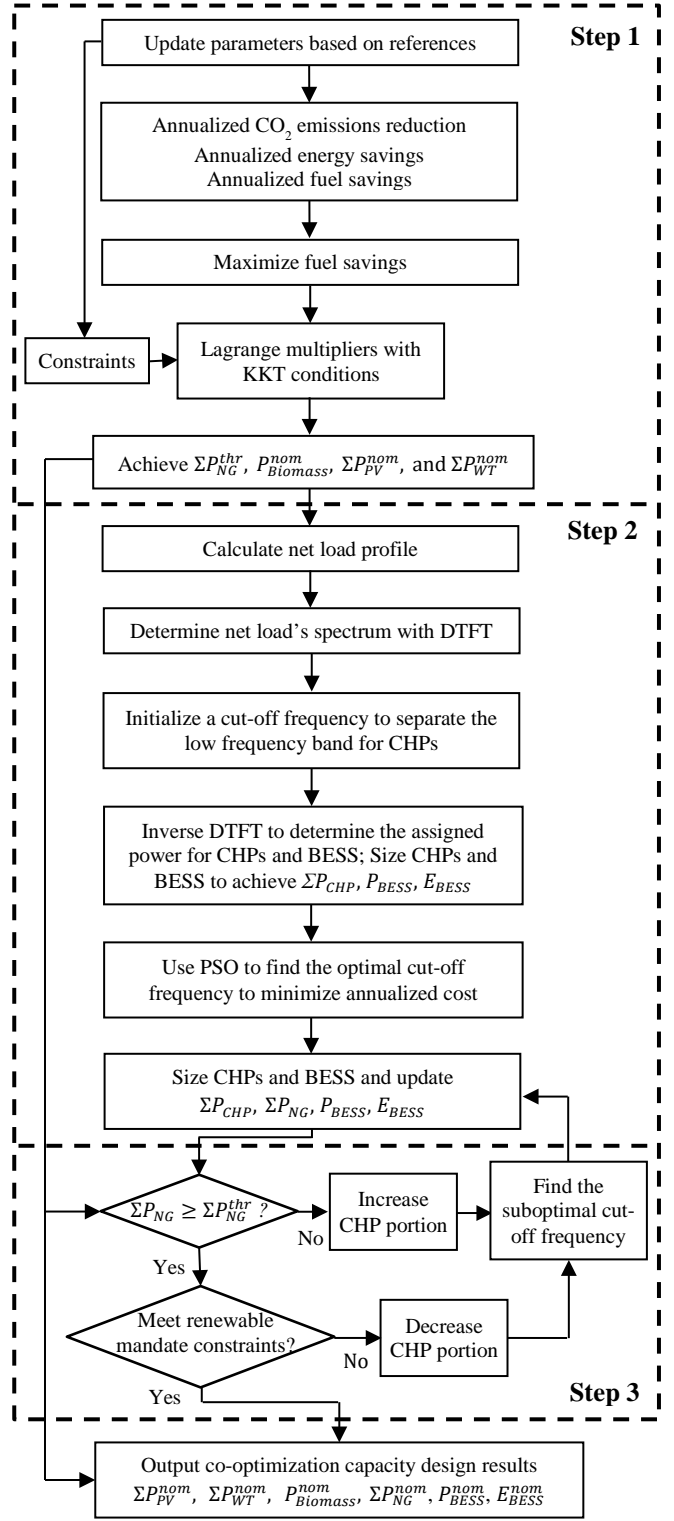

Fig. 2. Flowchart of co-optimization scheme for energy resource planning

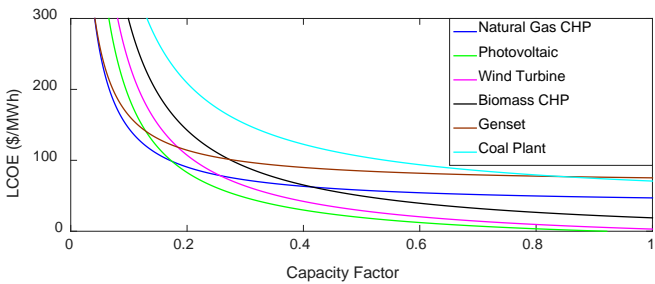

Fig. 1. LCOE vs. capacity factor curves for various types of DERs

### A. Step 1: Sizing of Renewable Energy Resources

According to the U.S. DOE requirement, community microgrids need to support the $CO_2$ emissions reduction and increase system energy efficiency [3], [4]. Many states also established renewable energy goals [5]. The design approach in this part uses Lagrange multipliers method with Karush-Kuhn-Tucker (KKT) conditions. In this step, the maximization of fuel savings serves as an objective function for sizing of renewable resources. This is because the intent is to exceed meeting the renewable portfolio standards by pursuing higher utilization of renewables—because of their lower LCOE and being environmental friendly. Besides, the threshold power rating of natural gas CHPs, $P_{NG}^{thr}$, is also achieved.

Before proceeding further, it should be noted that the community is assumed to have a wastewater treatment plant with an at most 0.5MW biomass CHP installed, $P_{Biomass}^{nom}$. This is because biomass is produced from the process of wastewater disposal, and is typically limited to 0.5 MW capacity. In addition, the capacity factor—the ratio of the actual average power output of an energy resource to its rated peak power if it could be run at full capacity—is accounted in calculating the annual fuel savings as well as the constraints considered in step 1 [27]. Depending on the weather conditions, the capacity factors of PV and WT are much smaller than other energy resources. However, as the fossil fuels are not consumed, they have low (even zero) $CO_2$ emissions. In the meanwhile, if considering renewable electricity production tax credit, their LCOEs are very competitive.

#### 1) Objective Function—Fuel Savings: 
Fuel consumption affects system operation cost, which varies as the operation schedule changes. The annual fuel savings (FS) is calculated as

$$FS = \gamma_{PV} \cdot \Sigma P_{PV}^{nom} + \gamma_{WT} \cdot \Sigma P_{WT}^{nom} + \gamma_{Biomass} \cdot P_{Biomass}^{nom} + \gamma_{NG} \cdot \Sigma P_{NG}^{thr} \quad (1)$$

where PV, WT, Biomass, and NG denote photovoltaic, wind turbine, biomass CHP, and natural gas CHP, respectively.

#### 2) $CO_2$ Emissions Reduction: 
The U.S. DOE requirement of $CO_2$ emissions reduction (ER) is given by

$$ER = \alpha_{PV} \cdot \Sigma P_{PV}^{nom} + \alpha_{WT} \cdot \Sigma P_{WT}^{nom} + \alpha_{NG} \cdot \Sigma P_{NG}^{thr} + \alpha_{Biomass} \cdot P_{Biomass}^{nom} \quad (2)$$

#### 3) System Energy Savings: 
The annual system energy savings (SES) is calculated as

$$SES = \beta_{PV} \cdot \Sigma P_{PV}^{nom} + \beta_{WT} \cdot \Sigma P_{WT}^{nom} + (\beta_e + \beta_{th}) \cdot (\Sigma P_{NG}^{thr} + P_{Biomass}^{nom}) \quad (3)$$

#### 4) Constraints: 
Based on the assumption of biomass CHP near wastewater treatment plant,

$$P_{Biomass}^{nom} \leq 0.5 \; MW \quad (4)$$

As per the requirements of U.S. DOE and various states, the constraints on renewables (5) and (6), $CO_2$ emissions (7), and system energy efficiency (8) are developed.

$$\frac{\Sigma P_{PV}^{nom}}{\Sigma P_{PV}^{nom} + \Sigma P_{WT}^{nom}} \geq \theta \; or \; \frac{\Sigma P_{PV}^{nom}}{\Sigma P_{PV}^{nom} + \Sigma P_{WT}^{nom} + P_{Biomass}^{nom}} \geq \theta \quad (5)$$

$$\frac{\Sigma P_{PV}^{nom} + \Sigma P_{WT}^{nom} + P_{Biomass}^{nom}}{\Sigma P_{PV}^{nom} + \Sigma P_{WT}^{nom} + \Sigma P_{NG}^{nom} + P_{Biomass}^{nom}} \geq \rho \quad (6)$$

$$\frac{ER}{E_{CO_2} \cdot L} \geq \omega \quad (7)$$

$$\frac{SES}{E_l + E_{th}} \geq \delta \quad (8)$$

where $E_l = \int_0^{8760} P_l(t)dt$ ; $E_{th} = \int_0^{8760} P_{th}(t)dt$

#### 5) Lagrange Multipliers with KKT Conditions: 
The Lagrange multipliers theory (LMT) is used for optimally sizing the renewable resources to achieve the maximal fuel savings. A general problem is shown in (9) and corresponding LMT formulation is presented in (10). The optimal solution $x^*$ is the local minimum/maximum for $f(x)$, achieved after satisfying the necessary and sufficient conditions of KKT in (11) and (12).

$$Minimize/Maximize \; f(x)$$
$$s.t. \; h_i(x) = 0, i = 1,2,\cdots,m \quad (9)$$
$$g_j(x) \leq 0, j = 1,2,\cdots,n$$

$$L(x,\lambda,\mu) = f(x) + \sum_{i=1}^{m} \lambda_i h_i(x) + \sum_{j=1}^{n} \mu_j g_j(x) \quad (10)$$

$$\nabla_x L(x^*,\lambda^*,\mu^*) = \nabla f(x^*) + \sum_{i=1}^{m} \lambda_i^* \nabla h_i(x^*) + \sum_{j=1}^{n} \mu_j^* \nabla g_j(x^*) = 0 \quad (11)$$

$$y^T \nabla_{xx}^2 L(x^*,\lambda^*,\mu^*)y \geq 0$$
$$\forall y \in R \; s.t. \; \nabla h_i(x^*)^T y = 0, \nabla g_j(x^*)^T y = 0 \quad (12)$$

### B. Step 2: Sizing of CHPs and BESS

To make full use of renewable energy and save operational costs, the renewable energy resources are always allowed to output maximum available power. However, because of their intermittent nature, the renewable power generation is uncertain and fluctuating. In addition, the load demand is often non-deterministic. Mathematical models of load and renewables are employed to mimic their stochastic behaviors. In this work, it is assumed that both PV and WTs operate at the unity power factor. Dispatchable generation units like CHPs and BESS should be capable of meeting the net load, equal to the power mismatch between renewable supply and load demand. A DTFT and PSO based coordinated sizing scheme is employed to assign the power allocation for minimizing annualized cost. Furthermore, a suitable reserve margin is necessary for satisfactory microgrid operation.

*1) Stochastic Model of Load:* Load demand is a stochastic variable. It is assumed to follow normal distribution [28].

$$f_{P_l}(P_l \mid \mu, \sigma^2) = \frac{1}{\sqrt{2 \cdot \sigma^2 \cdot \pi}} \cdot e^{-(P_l - \mu)^2 / 2 \cdot \sigma^2} \quad (13)$$

*2) Stochastic Model of PV:* According to statistical data, in each time frame, the solar irradiance is assumed to follow a beta distribution given by the following probability density function (PDF) [9], [29]:

$$f_{ir}(ir; a, b) = \frac{\Gamma(a+b)}{\Gamma(a) \cdot \Gamma(b)} \cdot \left(\frac{ir}{ir_{max}}\right)^{a-1} \cdot \left(1 - \frac{ir}{ir_{max}}\right)^{b-1} \quad (14)$$

where,
$\Gamma(a+b) = (a+b-1)!$
$\Gamma(a) = (a-1)!$
$\Gamma(b) = (b-1)!$

Since the total active power $P_{pv} = ir \cdot Area_{PV} \cdot \eta_{PV}$, the PDF of power available from PV can be deduced as follows:

$$f_{P_{PV}}(P_{PV}; a, b) = \frac{\Gamma(a+b)}{\Gamma(a) \cdot \Gamma(b)} \cdot \left(\frac{P_{PV}}{P_{PV\_max}}\right)^{a-1} \cdot \left(1 - \frac{P_{PV}}{P_{PV\_max}}\right)^{b-1} \quad (15)$$

*3) Stochastic Model of Wind:* Similar to solar irradiance, wind speed is stochastic and assumed to follow Weibull distribution [9], [29].

$$f_v(v; \tau, K) = \frac{K}{\tau} \cdot \left(\frac{v}{\tau}\right)^{K-1} \cdot e^{-\left(\frac{v}{\tau}\right)^K}, \quad v \geq 0 \quad (16)$$

$$F_v(v; \tau, K) = 1 - e^{-\left(\frac{v}{\tau}\right)^K} \quad (17)$$

A typical wind turbine starts to generate power when wind speed is larger than the cut-in speed, $V_{in}$, and the power output linearly increases with wind speed from the cut-in speed, $V_{in}$, to rated speed, $V_R$. When wind speed is between the rated value and the cut-out speed, $V_{out}$, the power output is the rated power, $P_R$. Otherwise, it does not generate power.

$$P_{wind}(v) = \begin{cases} P_R \cdot \dfrac{v - V_{in}}{V_R - V_{in}} & V_{in} \leq v \leq V_R \\ P_R & V_R \leq v \leq V_{out} \\ 0 & else \end{cases} \quad (18)$$

Thus the PDF of power available from wind turbine is presented below in (19).

*4) Net Load Profile:* The aim of this work is to size DERs in community microgrids. Therefore, the net load profile is generated with historical data based stochastic models. With the help of such data, the power capacity planning for DERs would cover prevailing year-round conditions. By selecting representative weekdays and weekends in four seasons, the scheme is able to save on computation time without compromising accuracy. However, the uncertainties of renewable energy and load demand lead to the indeterminacy of net load and make historical information inadequate. Therefore, this paper employs stochastic models to mimic net load behaviors and improve the performance of DERs sizing. The shaping and scaling of stochastic model parameters can be calculated with historical data. The net load profile can be expressed as load demand, $\Sigma P_l(t)$, minus renewable power output, $\Sigma P_{re}(t)$.

$$P_{nl}(t) = \Sigma P_l(t) - \Sigma P_{re}(t) \quad (20)$$

*5) Preliminary Sizing of CHPs and BESS with Discrete-Time Fourier Transform (DTFT):* The net load can be divided into two parts: a) Component with large power that varies smoothly over a longer duration, and b) small but frequently fluctuating power component. CHPs could take the smooth (i.e., flat) power variation and the BESS can compensate the small and frequent changes. Therefore, the aim of this step is to divide net load power between CHPs and BESS.

$$P_{nl}(t) = \Sigma P_{CHP}(t) + P_{BESS}(t) \quad (21)$$

As illustrated in Fig. 2, the DTFT is used to get the net load's spectrum. Then the cut-off frequency is initialized in a random way to assign the low-frequency part of the net load to CHPs and let BESS take care of high-frequency power fluctuations. This helps in lowering capital costs, as the BESS capacity, in terms of power and energy, will not be oversized. Once the allocation for CHPs is achieved, based on the cut-off

$$f_{P_{wind}}(P_{wind}; \tau, K) = \begin{cases} F_1 & P_{wind} = 0 \\ \dfrac{V_R - V_{in}}{P_R} \cdot \dfrac{K}{\tau} \cdot \left(\dfrac{V_{in} + (V_R - V_{in}) \cdot \frac{P_{wind}}{P_R}}{\tau}\right)^{K-1} \cdot e^{-\left(\frac{V_{in} + (V_R - V_{in})\frac{P_{wind}}{P_R}}{\tau}\right)^K} & 0 < P_{wind} < P_R \\ F_2 & P_{wind} = P_R \end{cases} \quad (19)$$

where,
$F_1 = 1 - [F_v(V_{out}; \tau, K) - F_v(V_{in}; \tau, K)]$
$F_2 = F_v(V_{out}; \tau, K) - F_v(V_R; \tau, K)$

frequency, the inverse DTFT is used to set the power share for CHPs in the time domain. However, the inverse DTFT can give negative values. Since the CHPs are incapable of absorbing power, it is necessary to shift the negative power to the BESS.

$$\Sigma P_{CHP}[k] = \begin{cases} \Sigma P_{CHP}[k], & \Sigma P_{CHP}[k] \geq 0 \\ 0, & \Sigma P_{CHP}[k] < 0 \end{cases} \quad (22)$$

After making such an adjustment, the BESS power share can be determined as the net load minus the CHPs assigned power. The final procedure is to size CHPs and BESS. As shown in (23), the power capacity of CHPs needs to meet the maximum power demand and also include a reserve margin to ride through forecast errors and unexpected events. In this paper, BESS efficiencies during both charging and discharging are assumed equal. So the BESS power capacity is calculated according to (24) and (25). It is discharging when $P_{BESS}[k] > 0$ and charging when $P_{BESS}[k] < 0$. Furthermore, the BESS energy capacity is determined in (26)-(27). In particular, the difference in stored energy from the original status to $k^{th}$ sample point is given by (26). As the BESS undergoes equal charging and discharging cycles, it is essential that (26) is zero at the end of each day. Equation (27) presents the BESS energy capacity sizing, where $SOC_{max}$ and $SOC_{min}$ are two predetermined parameters.

$$\Sigma P_{CHP}^{nom} = max\{\Sigma P_{CHP}[k]\} \cdot (1 + R) \quad (23)$$

$$P_{BESS}[k] = \begin{cases} P_{BESS}[k]/\eta_{BESS}, & P_{BESS}[k] \geq 0 \\ P_{BESS}[k] \cdot \eta_{BESS}, & P_{BESS}[k] < 0 \end{cases} \quad (24)$$

$$P_{BESS}^{nom} = max\{|P_{BESS}[k]|\} \quad (25)$$

$$E[k] = \sum_{i=0}^{k} (P_{BESS}[i] \cdot T), \qquad k = 0, \dots, N_S \quad (26)$$

$$E_{BESS}^{nom} = \frac{max\{E[k]\} - min\{E[k]\}}{SOC_{max} - SOC_{min}} \quad (27)$$

*6) Optimal Cut-off Frequency with PSO:* After their preliminary sizing, the power allocations of CHPs and BESS are estimated. However, since the initial cut-off frequency was chosen randomly, this may not be an optimal power assignment. Therefore, the particle swarm optimization (PSO) is used for determining the cut-off frequency to achieve the minimum annualized cost. The total annualized cost shown in (28) includes all the DER costs. For each type of DER, these costs include capital investment, operation, and maintenance (O&M) cost, fuel cost and tax credit, as presented in (29). The brackets in (29) indicate that fuel cost only applies to fuel based DERs and renewables have a tax credit as payback. The fuel cost of natural gas CHPs is calculated by (30). In this study, it is assumed that all natural gas CHPs share the load evenly; and therefore, their operational efficiencies are the same. However, this can be changed whenever necessary. For renewable electricity generation, the tax credit payback is expressed in (31). In the end, each DER's annualized capital cost is calculated using (32), considering lifetime and discount rate. The discount rate refers to the interest rate used in discounted cash flow (DCF) analysis to set the present value of future cash flows.

$$C_{total} = C_{PV} + C_{WT} + C_{Biomass} + C_{NG} + C_{BESS} \quad (28)$$

$$C = C^{capital} + C^{O\&M} (+ C^{fuel})(- C^{TaxCredit}) \quad (29)$$

$$C_{NG}^{fuel} = F \cdot \sum_{NG} \left( \frac{\sum_i P_{NG}[i] \cdot T}{\eta_{NG}[i]} \right) \quad (30)$$

$$C_{re}^{TaxCredit} = TaxCredit \cdot \sum_{re} \left( \sum_i P_{re}[i] \cdot T \right) \quad (31)$$

$$C^{capital} = Total\ Capital\ Cost \cdot \frac{r \cdot (1 + r)^y}{(1 + r)^y - 1} \quad (32)$$

### C. Step 3: Parity Check

After determining the power ratings of natural gas CHPs and BESS in step 2, a parity check is performed to see whether the constraints in step 1 are still satisfied. For this purpose, the calculated natural gas CHPs power capacity, $\Sigma P_{NG}$, is compared with its threshold power rating, $\Sigma P_{NG}^{thr}$. When $\Sigma P_{NG} < \Sigma P_{NG}^{thr}$, the constraints of $CO_2$ emissions and system energy efficiency are violated. Under this condition, the suboptimal cut-off frequency is chosen toward increasing CHP capacity, and the BESS and CHPs are resized. In addition, when $\Sigma P_{NG} \geq \Sigma P_{NG}^{thr}$, the state renewable mandate constraints may not be fully met. In such a case, the suboptimal cut-off frequency is selected toward reducing CHP capacity, and the BESS and CHPs are resized. This process is repeated until both $\Sigma P_{NG} \geq \Sigma P_{NG}^{thr}$ and the renewable mandate constraints are satisfied.

### D. Evaluation Indices

An evaluation system is necessary to weigh the benefits of a microgrid system. For this purpose, the following technical and economic indices are proposed.

*1) Renewable Energy Utilization:* With the higher penetration of renewable energy resources, the fuel consumption and $CO_2$ emissions can be reduced. The index of renewable energy utilization is defined as

$$URe = \frac{\int_0^{8760} (\Sigma P_{PV}(t) + \Sigma P_{WT}(t) + P_{Biomass}(t))dt}{\int_0^{8760} \Sigma P_l(t)dt}$$
$$or \quad URe = 1 - \frac{\int_0^{8760} \Sigma P_{NG}(t)dt}{\int_0^{8760} \Sigma P_l(t)dt} \quad (33)$$

*2) Fuel Savings:* In this paper, the objective function, in step 1, is the maximization of fuel savings with renewable resources (cf. Fig. 2). The expression for annual fuel savings is presented in (1).

*3) Regulations/Requirements:* The community microgrid has to meet the government regulations and customer requirements. Specifically, the DOE requirements and state renewable mandates, given in (5)-(8), are considered as the constraints of the co-optimization problem.

*4) Annualized Cost Evaluation:* The total annualized cost is minimized in step 2 by optimally sizing CHPs and BESS (cf.

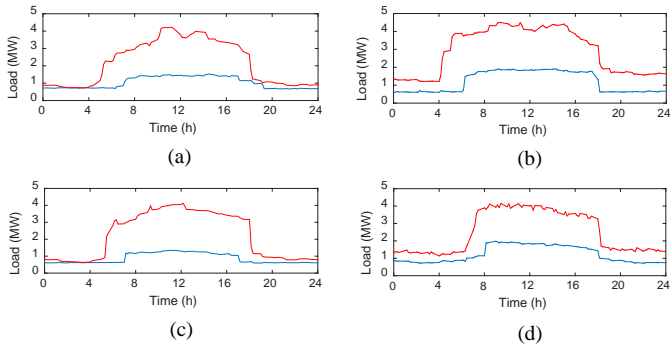

Fig. 3 Load profiles in four seasons: (a) spring, (b) summer, (c) autumn, (d) winter. (Note: red line is the load profile on weekdays, blue line is the load profile on weekends)

Fig. 2). The expressions for annualized cost were presented in (28)-(32).

## IV. CASE STUDY AND PERFORMANCE ASSESSMENT

In this section, the proposed co-optimization strategy for sizing of distributed energy resources is carried out for supporting the critical loads of a community in the State of Ohio, USA.

### A. Community Microgrid

The community chosen is a village with a 4.5 MW peak critical load [30]. The critical load is a commercial establishment, whose load data at 15-minute sample interval was provided by the local utility, AEP Ohio. Fig. 3 shows the representative load profiles of the commercial load on weekdays and weekends in four seasons. In the community microgrid, the reserve margin of CHPs is 10%, BESS system efficiency is 85% and the SOC limit is set as [50%, 100%]. Besides, the discount rate is 5% and the overall system lifetime is assumed as 20 years.

As observed in Fig. 3, the critical load profile of weekdays is much larger than that corresponding to weekends. Therefore, the generation capacity designed for weekdays is adequate for weekends. From [31], [32], the seasonal representative data of global horizontal irradiance (GHI) and wind speed could be selected to form stochastic models of load, PV, and WT.

Table I shows the coefficients of $CO_2$ emissions reduction, system energy savings and fuel savings for Ohio [33]–[37]. The base amount of annual $CO_2$ emissions in Ohio is 6330.55 tons/MW [37]. It should be noted that for CHPs the coefficients of annualized system energy savings contain electrical and thermal items, $\beta_e$ and $\beta_{th}$. Table II shows the capital costs, O&M costs and fuel costs of different DERs [38]–[40]. Besides, according to U.S. Internal Revenue Service, the renewable electricity production tax credit for PV and WT is 23 $/MWh and biomass's tax credit is 12 $/MWh

[26]. The annualized capital cost and O&M cost of sodium-sulfur (NaS) BESS, shown in Table III, are obtained from ES-Select™ tool [41] and references [42], [43]. As per the U.S. DOE requirements, the community microgrid should cut $CO_2$ emissions by ⩾20%, and increase system energy efficiency by ⩾20% [3], [4]. According to state renewable mandate for Ohio, the portion of renewable energy resources should be ⩾ 12.5% and solar (PV) has to be at least 0.5% [5].

### B. Coordinated Sizing of DERs

The proposed co-optimization scheme is tested using GAMS and MATLAB software. At first, the fuel savings maximization problem is formulated as a function of Lagrange multipliers with equality and inequality constraints. Once the necessary and sufficient KKT conditions are satisfied, the solution is guaranteed to be the optimal one. Then the nominal power rating of renewables and the threshold power rating for natural gas CHPs are determined with GAMS. Following the steps in the algorithm (cf. Fig. 2), implemented in MATLAB, DTFT is used to achieve the spectrum of net load power in the representative scenarios. Fig. 4 displays the critical load spectrum during weekdays in summer as an example. For the 15-minute sample interval, the frequency range of the spectrum is [0, 0.56 mHz]. The cut-off frequency divides the spectrum into two parts. The left portion (i.e., low frequency) is assigned to the CHPs and the right portion (i.e., high frequency) to the BESS. For this step, the particle swarm optimization (PSO) is used to find the optimal cut-off frequency within the frequency range—to derive the minimum annualized cost. Table IV displays results in four seasons, and the optimal cut-off frequencies are shown in Fig. 5. The red dots in Fig. 5 mark the lowest total annualized costs. Besides the optimal cut-off frequencies, the assigned capacities for CHPs and BESS are also presented in Table IV. Based on the observations and design considerations, the total power capacity of CHPs is determined as 2.04 MW (including 0.5 MW biomass CHP), and the BESS is sized at 1.24 MW and 3.25 MWh. This outcome is illustrated in Fig. 6.

Table II Cost parameters of DERs

| Cost<br>DER | Capital Cost ($/MW) | | O&M Cost ($/MW-year) | Fuel Cost ($/MWh) |
|---|---|---|---|---|
| | Total | Annualized | | |
| PV Panel | 2,800,000 | 224,700 | 19,000 | 0 |
| Wind Turbine | 1,710,000 | 137,200 | 29,000 | 0 |
| Natural Gas CHP | 1,000,000 | 80,000 | 91,000 | 10 |
| Biomass CHP | 1,000,000 | 80,000 | 91,000 | 0 |

Table III Cost parameters of BESS

| BESS | Annualized Capital Cost—Power ($/MW) | Annualized Capital Cost—Energy ($/MWh) | O&M Cost ($/MW/year) |
|---|---|---|---|
| NaS | 280,000 | 24,000 | 3000 |

Table I Coefficients update for the State of Ohio, USA

| DER<br>Coefficient | Renewables | | | Non-Renewables |
|---|---|---|---|---|
| | PV | WT | Biomass CHP | Natural Gas CHP |
| α (tons/MW) | 1479.7 | 1967.6 | 6437 | 6345 |
| β (MWh/MW) | 1594.3 | 2119.9 | 7008+9877.66<br>($\beta_e$)   ($\beta_{th}$) | 7008+9877.66<br>($\beta_e$)   ($\beta_{th}$) |
| γ (MMBtu/MW) | 23458.1 | 31191.7 | 103114 | 23340 |

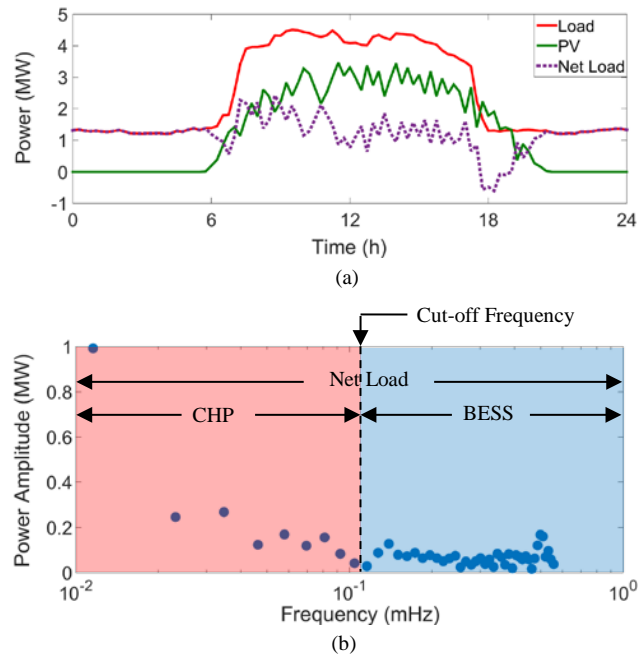

Fig. 4. Net load profile in (a) time domain; (b) frequency spectrum during a weekday of summer

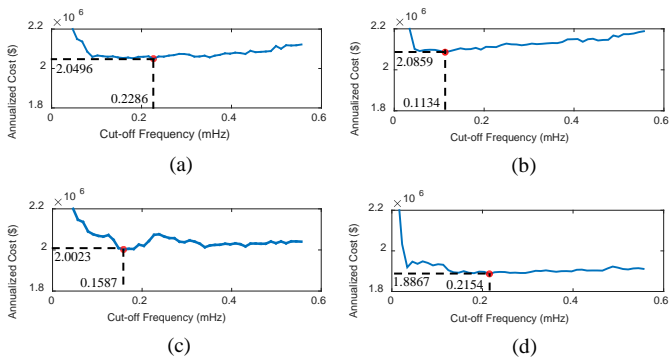

Fig. 5. Annualized cost versus cut-off frequency for weekdays in four seasons: (a) spring; (b) summer; (c) autumn; (d) winter.

### C. Performance Assessment

A performance assessment was carried out on the results

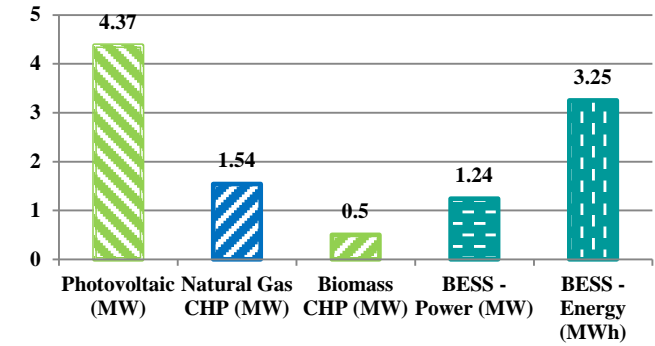

(a)

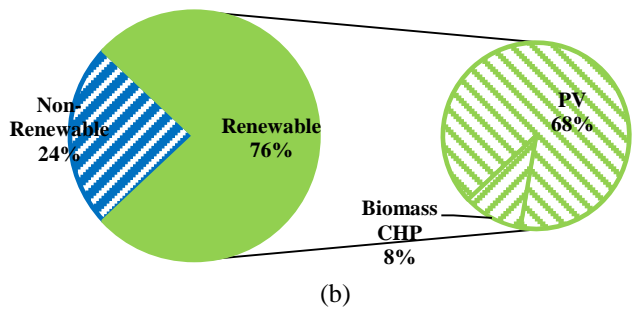

(b)

Fig. 6. Sizing of energy resources by the proposed scheme: (a) capacity of each resource; (b) portions of renewable and nonrenewable energy resources.

Table IV Coordinated sizing results of CHP and BESS in four seasons

| | Optimal Cut-off Frequency (mHz) | ΣCHP (MW) | BESS Power Capacity (MW) | BESS Energy Capacity (MWh) |
|---|---|---|---|---|
| Spring | 0.2268 | 1.97 | 1.06 | 1.72 |
| Summer | 0.1134 | 1.88 | 1.24 | 3.25 |
| Autumn | 0.1587 | 2.04 | 0.89 | 0.81 |
| Winter | 0.2154 | 1.86 | 0.68 | 1.23 |

Table V System evaluation

| System Evaluation Indices | Design Scheme / Proposed Scheme |
|---|---|
| Fuel Savings in MMBtu | 190,012.00 |
| System Energy Efficiency Increase | 27.20% |
| CO₂ Emissions Reduction | 26.70% |
| Utilization of Renewable Energy ($URe$) | 52.60% |
| DOE Requirements | ✔ |
| State Renewable Mandate | ✔ |
| Annualized Cost Minimization | ✔ |

Note: ✔ indicates satisfying the corresponding requirement

obtained for the case study. Fig. 6 presents the sizing results of distributed energy resources. Table V displays the system evaluation indices. Due to their zero-fuel-consumption nature and competitive LCOEs, renewables are sized to maximize fuel savings. It is also clear that the increases in system energy efficiency and CO₂ emissions reduction satisfy the DOE requirements. The annualized costs in cost ($) per MW load are displayed in Table VI. With high penetration of renewable energy resources, the capital cost takes over two-thirds of the total cost, while the fuel cost and O&M cost are very low.

For economic justification, the results from the co-optimization scheme are compared with two baseline cases in Fig. 7. Baseline I has only natural gas CHPs to just meet critical load demand and DOE requirements, and Baseline II includes the renewables and energy storage to also meet the state renewable energy mandate. When natural gas CHPs alone are used (i.e., Baseline I), the capital cost is much lower than the other two; however, its fuel cost and O&M cost are higher, taking up over 80% of the total cost. With the addition of renewables and energy storage (i.e., Baseline II), the capital cost is a large part, but its fuel cost and O&M cost are much lower. Although the total cost is higher in Baseline II than that in Baseline I, the state renewable energy mandate is met in Baseline II. As compared to Baseline II, the co-optimization scheme gives a better solution with lower total cost. This is because the co-optimization scheme maximizes fuel savings first to determine sizes of renewable energy resources and minimizes the annualized cost while satisfying the DOE requirements and state renewable energy mandate. Therefore, the capital cost of the co-optimization scheme is higher than Baseline I but lower than Baseline II. Although the fuel cost and O&M cost are a little higher than that in Baseline II, it is much lower than Baseline I's. Furthermore, the total annualized cost is lower than both Baselines I and II.

Table VI Annualized costs of the proposed scheme

| Annualized Cost | Design Scheme / Proposed Scheme |
|---|---|
| Capital Cost ($/MW-year) | 348,964.20 |
| O&M Cost ($/MW-year) | 60,531.10 |
| Fuel Cost ($/MW-year) | 65,344.40 |
| Tax Credit Payback ($/MW-year) | 47,290.03 |
| Total Cost ($/MW-year) | 427,549.67 |

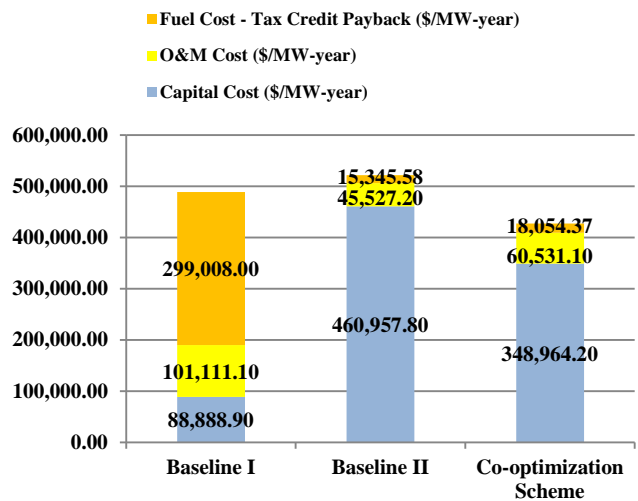

Fig. 7. Economic justification of the co-optimization scheme

Table VII DERs and BESS sizing comparison

| Results/Indices | Design Scheme Proposed Scheme | HOMER Pro |
|---|---|---|
| Natural Gas CHPs (MW) | 4.31 | 3.85 |
| Biomass CHP (MW) | 0.50 | 0.50 |
| BESS Power Capacity (MW) | 0.50 | 0.75 |
| BESS Energy Capacity (MWh) | 0.44 | 0.25 |

## V. COMPARISON WITH HOMER PRO

The results of proposed scheme are compared with HOMER Pro in this section. For HOMER Pro, the cost minimization is the sole objective, since the constraints of government regulations cannot be set. In other words, the results in Section IV cannot be directly used for the comparison. For a fair comparison, the proposed scheme is rerun for minimizing cost alone while ignoring the constraints on government regulations and fuel savings maximization, and the results are presented in Tables VII and VIII. As seen in Table VIII, the total cost of the proposed scheme is close but below that obtained using HOMER Pro. Furthermore, the proposed scheme considers the coordination between BESS and CHPs, thus helping better assign the net load demand. Although the CHPs have a slightly higher total capacity, they work more efficiently with BESS, and thereby lower the fuel cost. This is because the proposed scheme has higher energy capacity in BESS, and therefore it supports a longer operation of BESS to better smooth output from CHPs.

Comparing results in Table VI and Table VIII, with higher renewable penetration, the capital cost increases a lot while the fuel costs and O&M costs are reduced significantly.

Table VIII Annualized costs comparison assuming same conditions

| Annualized Cost | Design Scheme Proposed Scheme | HOMER Pro |
|---|---|---|
| Capital Cost ($/MW-year) | 118,968.90 | 125,333.30 |
| O&M Cost ($/MW-year) | 97,602.20 | 88,466.70 |
| Fuel Cost ($/MW-year) | 212,016.70 | 217,312.70 |
| Tax Credit Payback ($/MW-year) | 11,680.00 | 11,680.00 |
| Total Cost ($/MW-year) | 416,907.80 | 419,432.70 |

Nevertheless, the cost of renewable energy has decreased drastically over last decade, and this trend is expected to continue due to various technological developments. The proposed sizing methodology is a major step in this direction by integrating distributed energy resources, especially renewable energy resources, at the lowest annualized cost.

## VI. CONCLUSION

This paper presented a co-optimization scheme for distributed energy resource planning in community microgrids. The renewable energy resources are designed first to satisfy the U.S. DOE requirements and state renewable energy mandates. Besides, the maximization of fuel savings in the first step further aids in the utilization of renewable energy. Later, natural gas powered CHP units and BESS are sized in a coordinated approach to match the net load demand and to minimize total annualized cost. Evaluation indices including the renewables utilization were introduced. A case study was presented applying the proposed scheme for a community microgrid in the State of Ohio, USA. The impact of regulatory constraints on microgrid development was also examined. In addition, the proposed scheme was compared with HOMER Pro to confirm the results. The co-optimization scheme can be adapted well for different places and regulations.

## VIII. BIOGRAPHIES

**Chen Yuan** (S'13, M'17) received the B.S. degree in electrical engineering from Wuhan University, China, in 2012, and the Ph.D. degree in electrical and electronics engineering from The Ohio State University, Columbus, OH, USA, in 2016. He is a Postdoctoral Researcher at GEIRI North America. His current research interests include the high-performance computing, power systems analysis, control systems, and energy management systems.

**Mahesh S. Illindala** (S'01, M'06, SM'11) received the B.Tech. degree in electrical engineering from National Institute of Technology, Calicut, India, in 1995, the M.Sc.(Engg.) degree in electrical engineering from Indian Institute of Science, Bangalore, India, in 1999, and the Ph.D. degree in electrical engineering from the University of Wisconsin, Madison, WI, USA, in 2005.

From 2005 to 2011, Dr. Illindala was employed in Caterpillar R&D. Since 2011, he has been an Assistant Professor with the Department of Electrical and Computer Engineering, The Ohio State University, Columbus, OH, USA. Dr. Illindala is a recipient of the ONR Young Investigator Award in 2016. His research interests include microgrids, distributed energy resources, electrical energy conversion and storage, power system applications of multiagent systems, protective relaying and advanced electric drive transportation systems.

**Amrit Khalsa** (M'91) is a Staff Engineer at AEP's Dolan Technology Center. His research interests include microgrids, distributed energy resources, and smart grid technology from the distribution system down to the customer.

Mr. Khalsa is a registered Professional Engineer in the State of Ohio. He is the recipient of an EPRI Technology Transfer Award in 2012 for contributions to developing ANSI Standard CEA-2045, a modular communication interface for demand response.